\documentclass[12pt]{article}
\usepackage{amssymb,epic}
\evensidemargin\oddsidemargin
\advance\textheight-\headheight
\advance\textheight-\headsep
\advance\textheight-\topskip

\newtheorem{theorem}{Theorem}

\newcommand{\proof}{\noindent\underline{Proof.}\ }
\newtheorem{bem}[theorem]{Remark}

\newtheorem{proposition}[theorem]{Proposition}

\renewcommand{\setminus}{-}

\newcommand{\be}{\begin{enumerate}}
\newcommand{\ee}{\end{enumerate}}
\newcommand{\bi}{\begin{itemize}}
\newcommand{\ei}{\end{itemize}}

\newcommand{\ba}{\begin{array}}
\newcommand{\ea}{\end{array}}

\newcommand{\lsemi}{:}

\newcommand{\ev}{\mbox{\rm ev}}

\newcommand{\Z}{{\mathbb{Z}}}

\newcommand{\N}{{\mathbb{N}}}
\newcommand{\R}{{\mathbb{R}}}

\newcommand{\eb}{\phantom{zzz}\hfill{$\square $}\smallskip}

\renewcommand{\em}{\sf}

\title{Strongly modular lattices with long shadow.}

\author {
Gabriele Nebe \\
{\small Abteilung Reine Mathematik, Universit\"at Ulm, 89069 Ulm, Germany} \\
{\small nebe@mathematik.uni-ulm.de} 
}

\begin{document}

\maketitle

\small
%{\sc R\'esum\'e:}
%On donne la clssification des r\'eseaux fortement modulaires
%qui ont maximale ou 

{\sc Abstract:}\footnote{MSC 11H31, 11H50}
We classify strongly modular lattices with 
longest and second longest possible shadow.

\normalsize

\section{Introduction}

To an integral lattice $L $ in the euclidean space $(\R^n , (,))$, 
one associates the set of {\em characteristic vectors} $v \in \R^n$ 
with $(v,x) \equiv (x,x) \mbox{ mod } 2\Z $ for all $x\in L$.
They form a coset modulo $2 L^*$, where 
$$L^* = \{ v\in \R^n \mid (v,x) \in \Z \  \forall x\in L \} $$
is the {\em dual lattice} of $L$.
Recall that $L$ is called {\em integral}, if $L\subset L^*$ and 
{\em unimodular}, if $L=L^*$.
For a unimodular lattice, the square length of a characteristic vector
is congruent to $n$ modulo $8$ and there is always a
characteristic vector of square length $\leq n$.
In \cite{Elkies1} Elkies characterized 
the standard lattice $\Z^n$ as 
the unique unimodular lattice of dimension $n$, for which all 
characteristic vectors have square length $\geq n$. \cite{Elkies} 
gives the short list of unimodular lattices $L$ with 
$\min (L) \geq 2 $ such that all characteristic vectors of $L$ have 
length $\geq n-8$. The largest dimension $n$ is 23 and in dimension 23
this lattice is the {\em shorter Leech lattice} $O_{23}$ of minimum 3.
In this paper, these theorems are generalized to certain 
strongly modular lattices.
Following \cite{queb1} and \cite{queb2}, 
an integral lattice $L$ is called {\em $N$-modular}, if $L$ is isometric to
its rescaled dual lattice $\sqrt{N} L^*$.
A $N$-modular lattice $L$ is called {\em strongly $N$-modular},
if $L$ is isometric to all rescaled partial dual lattices
$\sqrt{m} L^{*,m} $, for all exact divisors $m$ of $N$, 
where $$L^{*,m} := L^* \cap \frac{1}{m} L .$$

The simplest strongly $N$-modular lattice is
$$C_{N} := \perp _{d\mid N} \sqrt{d} \Z $$
of dimension $\sigma_0(N) := \sum _{d \mid N} 1 $ the number of divisors 
of $N$.
The lattice $C_N$ plays the role of $\Z = C_1$ for square free $N>1$.

With the help of modular forms
Quebbemann \cite{queb2} shows that for 
$$N \in {\cal L} := \{ 1,2,3,5,6,7,11,14,15,23 \} $$
(which is the set of all positive integers $N$ such that 
the sum of divisors 
$$\sigma _1 (N) := \sum _{d \mid N } d $$
divides 24), 
the minimum of an even strongly $N$-modular lattice $L$ of dimension $n$
satisfies
$$\min (L ) \leq  2 + 2 \lfloor \frac{n \ \sigma_1(N) }{24 \  \sigma _0(N)} 
\rfloor .$$
Strongly modular lattices meeting this bound are called {\em extremal}.
Whereas Quebbemann restricts to even lattices, \cite{RaS} shows that the
same bound 
also holds for odd strongly modular lattices, where there is one exceptional
dimension $n =  \sigma_0(N) ( \frac{24 }{\sigma _1(N) } -1)$, where the bound on the minimum is 3 (and not 2).
In this dimension, there is a unique lattice $S^{(N)}$ of minimum 3.
For $N=1$, this is again the shorter Leech lattice $O_{23}$.
The main tool to get the bound for odd lattices is the
{\em shadow} 
$$S(L) := \{ \frac{v}{2} \mid v \mbox{ is a characteristic vector of } L \} .
 $$
If $L$ is even, then $S(L) = L^*$ and if $L$ is odd, $S(L) = L_0^* \setminus L^*$,
where $$L_0 := \{ v\in L \mid (v,v) \in 2\Z \} $$ is the even sublattice 
of $L$.

The main result of this paper is Theorem \ref{main}.
It is shown that for a strongly $N$-modular lattice $L$ that is rationally
equivalent to $C_N^k$, the minimum 
$$\mbox{$\min _0$} (S(L)) 
:= \min \{ (v,v) \mid v \in S(L) \}   $$ equals $$
M^{(N)}(m,k) := \left\{ \begin{array}{cc} 
 \frac{1}{N} (k\frac{\sigma_1(N)}{4} - 2m ) & \mbox{ if $N$ is odd } \\
 \frac{1}{N} (k\frac{\sigma_1(N/2)}{2} - m ) & \mbox{ if $N$ is even } \end{array} \right. $$
for some $m\in \Z _{\geq 0}$.
If $\min _0(S(L)) = M^{(N)}(0,k)$, then $L\cong C_N^k$.
For the next smaller possible minimum 
$\min _0(S(L)) = M^{(N)}(1,k)$ one gets that $L \cong C_N^l \perp L'$,
where $\min (L') > 1$ and $\dim (L') \leq \sigma_0(N) ( s(N) -1)$ for odd $N$
resp. 
 $\dim (L') \leq \sigma_0(N)  s(N) $ for even $N$.
The lattices $L'$ of maximal possible dimensions have minimum 3 and are uniquely determined:
 $L' = S^{(N)}$, if $N$ is odd and $L'=O^{(N)}$ (the ``odd analogue''
of the unique extremal strongly $N$-modular lattice of dimension $\sigma_0(N) s(N) $) if $N$ is even
(see \cite[Table 1]{RaS}).

The main tool to prove this theorem are the formulas for the
theta series of a strongly $N$-modular lattice $L$ and of its shadow 
$S(L)$  developed in \cite{RaS}. Therefore we briefly repeat these formulas in
the next section. 

\section{Theta series}

For a subset $S\subset \R ^n$, which is a finite union of cosets of an integral
 lattice
we put  its {\em theta series} 
$$\Theta _S(z) := \sum _{v\in S} q^{(v,v)}, \ \ q=\exp(\pi i z) .$$
The theta series of strongly $N$-modular lattices are modular forms
for a certain discrete subgroup $\Gamma _N$ of $SL_2(\R )$ (see \cite{RaS}).
Fix $N\in {\cal L}$ and put 
$$g_1^{(N)}(z):= \Theta _{C_N}(z) = 1+2q + 2 \ev(N)  q^2 + \ldots $$
where $$\ev(N) := \left\{ \begin{array}{ll} 1 & \mbox{ if $N$ is even } \\
0 & \mbox{ if $N$ is odd } \end{array} \right. .$$

Let $\eta $ be the Dedekind eta-function 
$$\eta (z) := q^{\frac{1}{12}} \prod _{m=1}^{\infty } (1-q^{2m}), \ \ q=\exp(\pi i z) .
$$
and put
$$\eta ^{(N)} (z) := \prod_{d\mid N} \eta (dz) .$$
If $N$ is odd define
$$g_2^{(N)}(z):= (\frac{\eta ^{(N)}(z/2) \eta ^{(N)}(2z) }{\eta ^{(N)}(z)^2} )^s$$
and if $N$ is even then
$$g_2^{(N)}(z):= (\frac{\eta ^{(N/2)}(z/2) \eta ^{(N/2)}(4z) }{\eta ^{(N/2)}(z)\eta ^{(N/2)}(2z)} )^s.$$
Then $g_2^{(N)}$ generates the field of modular functions of $\Gamma _N$.
It is a power series in $q$ starting with
$$g_2^{(N)}(z) = q - sq^2 + \ldots .$$

\begin{theorem}(\cite[Theorem 9, Corollary 3]{RaS})
{\label{RaS}}
Let $N\in {\cal L} $ and $L$ be a strongly $N$-modular lattice
that is rational equivalent to $C_N^k$.
Define $l_N:= \frac{1}{8} \sigma_1(N)$, if $N$ is odd and 
$l_N:= \frac{1}{6} \sigma_1(N)$, if $N$ is even.
Then 
$$\Theta _{L} (z) = g_1^{(N)}(z)^k \sum _{i=0}^{\lfloor k l_N \rfloor } c_i g_2^{(N)}(z)^i $$
for $c_i\in \R$. The theta series of the rescaled shadow $S:=\sqrt{N} S(L)$
of $L$ is
$$\Theta _{S} (z) = s_1^{(N)}(z)^k \sum _{i=0}^{\lfloor k l_N \rfloor } c_i s_2^{(N)}(z)^i $$
where $s_1^{(N)}$ and $s_2^{(N)}$
are the corresponding ``shadows'' of $g_1^{(N)}$ and $g_2^{(N)}$.
\end{theorem}

For odd $N$
 $$s_1^{(N)}(z) = 2^{\sigma_0 (N)} \frac{\eta ^{(N)}(2z)^2}{\eta ^{(N)} (z)} $$
and 
$$s_2^{(N)}(z) = -2^{-s(N)\sigma_0(N)/2} (\frac{\eta ^{(N)}(z)}{\eta ^{(N)}(2z)})^{s(N)} $$
For $N=2$ one has 
$$s_1^{(2)}(z) = \frac{2\eta (z)^5\eta (4z)^2}{\eta(z/2)^2\eta(2z)^3} $$
and
$$s_2^{(2)}(z) = -\frac{1}{16} (\frac{\eta (z/2)\eta (2z)^2}{\eta(z)^2\eta(4z)})^8 $$
which yields $s_1^{(N)}$ and $s_2^{(N)} $ for $N=6,14 $ as
$$s_1^{(N)} = s_1^{(2)}(z) s_1^{(2)}(\frac{N}{2} z) $$ and
$$s_2^{(N)} = (s_2^{(2)}(z) s_2^{(2)}(\frac{N}{2} z))^{s(N)/s(2)} .$$ 
If $N$ is odd, then $s_1^{(N)}$ starts with $q^{\sigma _1(N)/4}$ and
 $s_2^{(N)}$ starts with $q^{-2}$.
If $N$ is even, then $s_1^{(N)}$ starts with $q^{\sigma_1(\frac{N}{2})/2} $ and 
$s_2^{(N)}$ starts with $q^{-1}$.

\section{Strongly modular lattices with long shadow.}

\begin{proposition}{\label{abs}}
Let $N\in \N$ be square free and let $L$ be a strongly $N$-modular lattice.
If $L$ contains a vector of length 1, then $L$ has an orthogonal summand $C_{N}$.
\end{proposition}

\proof
Since $L$ is an integral lattice that contains a vector of length 1, the unimodular
lattice $\Z $ is an orthogonal summand of $L$.
Hence $L = \Z \perp L' $.
If $d$ is a divisor of $N$,
% then $d$ is an exact divisor of $N$, because $N$ is square free.
then $$L \cong \sqrt{d} L^{*,d} = \sqrt{d} \Z \perp \sqrt{d} (L')^{*,d} $$
by assumption. Hence $L$ contains an orthogonal summand $\sqrt{d}\Z $ for all
divisors $d$ of $N$ and
therefore $C_{N}$ is an orthogonal summand of $L$.
\eb

\begin{theorem} \label{main} (see \cite{Elkies} for $N=1$)
Let $N\in {\cal L} $ and  $L$ be a strongly $N$-modular lattice
that is rational equivalent to $C_N^k$.
Let $M^{(N)}(m,k)$ be as defined in the introduction.
\begin{itemize}
\item[(i)]
$\min _0(S(L) )  = M^{(N)}(m,k)$
for some $m\in \Z _{\geq 0}$.
\item[(ii)]
If 
$\min _0( S(L) )  = M^{(N)}(0,k) $ then 
$L \cong C_N^k$.
\item[(iii)]
If 
$\min _0( S(L) )  = M^{(N)}(m,k)$ then 
$L \cong C_N^a \perp L'$, where  $L'$ is a strongly $N$-modular lattice 
rational equivalent to $C_N^{k-a}$ with 
\\
$\min (L' ) \geq 2$ and 
$\min _0( S(L') )  = M^{(N)}(m,k-a)$.
\item[(iv)]
If 
$\min _0( S(L) )  = M^{(N)}(m,1)$ and $\min (L) \geq 2$,
then the number of vectors of length $2$ in $L$ is
$$2k(s(N)+\ev (N) - (k+1)) .$$
In particular $k \leq k_{max}(N)$ with $$k_{max}(N) = s(N) -1 +\ev (N) $$ 
and if $k=k_{max}(N)$, then $\min (L) \geq 3$.
\end{itemize}
\end{theorem}

\proof
(i) follows immediately from Theorem \ref{RaS}.
\\
(ii) In this case the theta series of $L$ is $g_1^k$. In particular $L$ contains
$2k$ vectors of norm 1. Applying Proposition \ref{abs} one finds that
$L\cong C_N$.
\\
(iii)
Follows from Proposition \ref{abs} and Theorem \ref{RaS}.
\\
(iv)
Since $\min (L) >1$, $\Theta _L = g_1^k - 2k g_1^kg_2 $.
Explicit calculations give the number of norm-2-vectors in $L$.
\eb

The following table gives the maximal dimension 
$n_{max} (N)= \sigma_0(N) k_{max}(N) $ of a lattice in Theorem \ref{main} (iv).
$$\begin{array}{|c|c|c|c|c|c|c|c|c|c|c|}
\hline
N & 1 & 2 & 3 & 5 & 6 & 7 & 11 & 14 & 15 & 23 \\
\hline
\sigma_1(N)  & 1 & 3 & 4 & 6 & 12 & 8 & 12 & 24 & 24 & 24 \\
\hline
k_{max}(N) & 23 & 8 & 5 & 3 & 2 & 2 & 1 & 1 & 0 & 0 \\
\hline
n_{max} (N)& 23 & 16 & 10 & 6 & 8 & 4 & 2 & 4 & 0 & 0 \\
\hline
\end{array}
$$

The lattices $L$ with $\min _0 ( S(L) ) = M(1,k)$ are listed in an appendix.
These are only finitely many since $k$ is bounded by $k_{max}$.
In general it is an open problem whether for all $m$, 
there are only finitely many strongly $N$-modular lattices $L$ rational equivalent 
to $C_N^k$ for some $k$ and of minimum $\min (L) > 1$ such that 
$\min _0 ( S(L) ) = M(m,k)$. 
For $N=1$, Gaulter \cite{Gaulter} proved that 
$k \leq 2907$ for $m=2$ and $k\leq 8388630$ for $m=3$.

\begin{theorem} (cf. \cite{Elkies} for $N=1$)
Let $N\in {\cal L}$ be odd and $k\in \N$ such that 
$$\frac{8}{\sigma_1(N)} \leq k\leq k_{max}(N) = \frac{24}{\sigma_1(N)} - 1 .$$
 Then there is a unique strongly $N$-modular lattice 
$L := L_k(N)$ that is rational equivalent to $C_N^k$ such that 
$\min (L) > 1 $ and $\min _0 ( S(L) ) = M^{(N)}(1,k) $,
except for $N=1$, where there is no such lattice in dimension 9, 10, 11, 13
and there are two lattices in dimension 18 and 20 (see \cite{Elkies}).
If $k=k_{max}(N)$, then $L$ is the shorter lattice $L=S^{(N)}$ described in
\cite[Table 1]{RaS} and $\min(L) = 3$.
\end{theorem}

\proof
For $N=15$ and $N=23$ there is nothing to show since $k_{max} (N) = 0$.
The case $N=1$ is already shown in \cite{Elkies}.
It remains to consider $N \in \{ 3,5,7,11 \} $.
Since $N$ is a prime, there are only 2 genera of strongly modular lattices, 
one consisting of even lattices and one of odd lattices.
With a short MAGMA program using Kneser's neighboring method, one obtains a list of 
all lattices in the relevant genus.
In all cases there is a unique lattice with the right number of vectors of length 2.
Gram matrices of these lattices are given in the appendix.
\eb

\begin{bem}
For $N=1$ and dimension $n=9,10,11$ the theta series of the 
hypothetical shadow has non integral resp.
odd coefficients, so there is no lattice $L_{n}(1)$.
\end{bem}

\begin{theorem}{\label{Neven}}
Let $N\in {\cal L}$ be even and $k\in \N$ such that
$$\frac{2}{\sigma_1(N/2)} \leq k\leq k_{max}(N) = \frac{24}{\sigma_1(N)} .$$
If $(k,N) \neq (3,2)$ then there are  strongly $N$-modular lattices
$L := L_k(N)$ that are rational equivalent to $C_N^k$ such that
$\min (L) > 1 $ and $\min _0 ( S(L) ) = M^{(N)}(1,k) $,
If $k=k_{max}(N)$, then $L_k(N)$ is unique.
It is the odd lattice $L=O^{(N)}$ described in
\cite[Table 1]{RaS}  and $\min(L) = 3$.
\end{theorem}

\begin{bem}
For $N=2$ and  $k=3$ the corresponding shadow modular form has non integral 
coefficients, so there is no lattice $L_{3}(2)$.
\end{bem}

\begin{bem}
All odd lattices $L_k(N)$  in Theorem \ref{Neven} lie in the genus
of $C_N^k$.
\end{bem}

\section{Appendix: The lattices $L_k(N)$.}

{\bf The lattices $L_k(1)$:}

The lattices $L_k(1)$ are already listed in \cite{Elkies}. 
They are uniquely determined by their root-sublattices $R_k$ 
and given
in the following table:
\begin{center}
\begin{tabular}{|c|c|c|c|c|c|c|c|c|c|c|c|c|}
\hline
$k$ & 8 & 12 & 14 & 15 & 16 & 17 & 18 & 19 & 20 & 21 & 22 & 23 \\ 
\hline
$R_k$ & $E_8$ &  $D_{12}$ & $E_7^2$ & $A_{15}$ & $D_8^2$ & $A_{11}E_6$ &
$D_6^3$, $A_9^2$ & $A_7^2D_5$ & $D_4^5$, $A_5^4$ & $A_3^7$ & $A_1^{22}$ & $0$ \\
\hline
\end{tabular}
\end{center}

{\bf The lattices $L_k(N)$ for $N>1$ odd:}
\begin{itemize}
\item[$L_2(3)$:]
$\left( \begin{array}{cc}  2 & 1 \\ 1 & 2 \end{array} \right) 
\perp \left( \begin{array}{cc}  2 & 1 \\ 1 & 2 \end{array} \right)  \cong A_2\perp A_2$,
Automorphism group: $D_{12} \wr C_2$.
\item[$L_3(3)$:]
\small
$\left(\begin{array}{@{}r@{}r@{}r@{}r@{}r@{}r@{}} 
2 & 1 & 1 & 1 & 1 & 1 \\
1 & 2 & 1 & 1 & 1 & 1 \\
1 & 1 & 2 & 1 & 1 & 1 \\
1 & 1 & 1 & 3 & 0 & 0 \\
1 & 1 & 1 & 0 & 3 & 0 \\
1 & 1 & 1 & 0 & 0 & 3 \end{array} \right) $ 
\normalsize
Automorphism group: order 1152.
\item[$L_4(3)$:]
\small
$\left(\begin{array}{@{}r@{}r@{}r@{}r@{}r@{}r@{}r@{}r@{}} 
2&0&0&0&\mbox{-}1&\mbox{-}1&0&1\\
0&2&0&0&1&1&1&0\\
0&0&2&0&0&\mbox{-}1&\mbox{-}1&\mbox{-}1\\
0&0&0&2&\mbox{-}1&0&\mbox{-}1&\mbox{-}1\\
\mbox{-}1&1&0&\mbox{-}1&3&1&1&0\\
\mbox{-}1&1&\mbox{-}1&0&1&3&1&0\\
0&1&\mbox{-}1&\mbox{-}1&1&1&3&1\\
1&0&\mbox{-}1&\mbox{-}1&0&0&1&3\end{array} \right) $
\normalsize
 Automorphism group: order 6144.
\item[$L_5(3)$:]
\small
$\left(\begin{array}{@{}r@{}r@{}r@{}r@{}r@{}r@{}r@{}r@{}r@{}r@{}} 
3 & \mbox{-}1 & \mbox{-}1 & \mbox{-}1 & 0 & \mbox{-}1 & \mbox{-}1 & \mbox{-}1 & 1 & 0\\
\mbox{-}1 & 3 & 1 & \mbox{-}1 & \mbox{-}1 & 1 & \mbox{-}1 & 1 & 1 & 0\\
\mbox{-}1 & 1 & 3 & 1 & 1 & 1 & \mbox{-}1 & 1 & 0 & \mbox{-}1\\
\mbox{-}1 & \mbox{-}1 & 1 & 3 & 1 & 0 & 1 & 1 & \mbox{-}1 & \mbox{-}1\\
0 & \mbox{-}1 & 1 & 1 & 3 & 1 & 0 & \mbox{-}1 & \mbox{-}1 & 0\\
\mbox{-}1 & 1 & 1 & 0 & 1 & 3 & \mbox{-}1 & 0 & \mbox{-}1 & \mbox{-}1\\
\mbox{-}1 & \mbox{-}1 & \mbox{-}1 & 1 & 0 & \mbox{-}1 & 3 & 0 & \mbox{-}1 & 1\\
\mbox{-}1 & 1 & 1 & 1 & \mbox{-}1 & 0 & 0 & 3 & 1 & 0\\
1 & 1 & 0 & \mbox{-}1 & \mbox{-}1 & \mbox{-}1 & \mbox{-}1 & 1 & 3 & 1\\
0 & 0 & \mbox{-}1 & \mbox{-}1 & 0 & \mbox{-}1 & 1 & 0 & 1 & 3\end{array} \right) $
\normalsize
Automorphism group: $\pm U_4(2) .2$ of order 103680
\item[$L_2(5)$:]
 $\left( \begin{array}{cc}  3 & 1 \\ 1 & 2 \end{array} \right) 
\perp \left( \begin{array}{cc}  3 & 1 \\ 1 & 2 \end{array} \right)  $,
Automorphism group: $(\pm C_2) \wr C_2$ of order 32.
\item[$L_3(5)$:]
\small
$\left(\begin{array}{@{}r@{}r@{}r@{}r@{}r@{}r@{}} 
3 & \mbox{-}1 & 1 & -1 & 1 & 0\\
\mbox{-}1 & 3 & \mbox{-}1 & 0 & 1 & 1\\
1 & \mbox{-}1 & 3 & 1 & 0 & 1\\
\mbox{-}1 & 0 & 1 & 3 & \mbox{-}1 & 1\\
1 & 1 & 0 & \mbox{-}1 & 3 & 1\\
0 & 1 & 1 & 1 & 1 & 3\end{array} \right) $.
\normalsize
 Automorphism group: $\pm S_5$ of order 240.
\item[$L_1(7)$:]
$\left(\begin{array}{cc} 
2 & 1 \\ 1 & 4 \end{array} \right) $. 
Automorphism group: 
$\pm C_2$.
\item[$L_2(7)$:]
\small
 $\left(\begin{array}{@{}r@{}r@{}r@{}r@{}} 
3 & \mbox{-}1 & 1 & 0\\
\mbox{-}1 & 3 & 0 & 1\\
1 & 0 & 3 & 1\\
0 & 1 & 1 & 3\end{array}  \right)$. 
\normalsize
Automorphism group: order 16.
\item[$L_1(11)$:]
 $\left(\begin{array}{cc} 
2 & 1 \\ 1 & 6 \end{array} \right) $. 
Automorphism group:
$\pm C_2$.
\end{itemize}

{\bf The lattices $L_k(N)$ for $N$ even:} \\
For $N=2$ there is only one genus of odd lattices to be considered.
Also for $N=14$ there is only one odd genus for each $k$, since 
2 is a square modulo 7.
For $N=6$, there are 2 such genera, since 
$L:=(\sqrt{2}\Z )^2 \perp (\sqrt{3}\Z )^2$ is not in the genus  of $C_6$. 
The genus of $L$ contains no strongly modular lattices.
The genus of $L\perp C_6 $  contains 3 lattices with minimum 3, none of which is
strongly modular.
\begin{itemize}
\item[$L_2(2):$] 
$L_2(2) = D_4 $ with automorphism group $W(F_4) $ of order 1152.
\item[$L_4(2):$] 
\small
$\left(\begin{array}{@{}r@{}r@{}r@{}r@{}r@{}r@{}r@{}r@{}} 
2 & 0 & 0 & 0 & 0 &\mbox{-}1 &\mbox{-}1&  1 \\
0 & 2 & 0 & 0 & 0 & 1 & 1& \mbox{-}1 \\
0 & 0 & 2 &\mbox{-}1 & 1 & 1 & 1& \mbox{-}1 \\
0 & 0 &\mbox{-}1 & 2 &\mbox{-}1 &\mbox{-}1 & 0&  1 \\
0 & 0 & 1 &\mbox{-}1 & 2 & 1 & 0& \mbox{-}1 \\
\mbox{-}1 & 1 & 1 &\mbox{-}1 & 1 & 3 & 2& \mbox{-}2 \\
\mbox{-}1 & 1 & 1 & 0 & 0 & 2 & 3& \mbox{-}1 \\
 1 &\mbox{-}1 &\mbox{-}1 & 1 &\mbox{-}1 &\mbox{-}2 &\mbox{-}1&  3 \end{array} \right) $.
\normalsize
$\begin{array}{lll} \mbox{
The root sublattice is $D_4 \perp A_1^4$} \\
\mbox{and the automorphism group of $L_4(2)$ is} \\
 \mbox{$W(F_4) \times (C_2 ^4 \lsemi D_8 ) $ of order 147456. } \end{array}
$
\item[$L_5(2)$:] 
\small
$\left(\begin{array}{@{}r@{}r@{}r@{}r@{}r@{}r@{}r@{}r@{}r@{}r@{}} 
2 & 1 & 1 & 1 & 1 & 1 & 1 & 1 & 1 & 1\\
1 & 2 & 1 & 1 & 1 & 1 & 1 & 1 & 1 & 1\\
1 & 1 & 2 & 1 & 1 & 1 & 1 & 1 & 1 & 1\\
1 & 1 & 1 & 2 & 1 & 0 & 0 & 0 & 0 & 0\\
1 & 1 & 1 & 1 & 2 & 1 & 1 & 1 & 1 & 1\\
1 & 1 & 1 & 0 & 1 & 3 & 1 & 1 & 1 & 1\\
1 & 1 & 1 & 0 & 1 & 1 & 3 & 1 & 1 & 1\\
1 & 1 & 1 & 0 & 1 & 1 & 1 & 3 & 1 & 1\\
1 & 1 & 1 & 0 & 1 & 1 & 1 & 1 & 3 & 1\\
1 & 1 & 1 & 0 & 1 & 1 & 1 & 1 & 1 & 3 \end{array} \right) $.
\normalsize
$\begin{array}{lll} \mbox{The root sublattice is $A_5$} \\ \mbox{and the 
automorphism group of $L_5(2)$ is } \\ \mbox{$\pm S_6 \times S_6$
of order 1036800.} \end{array} $
\item[$L_6(2):$]  
There are two such lattices:
\\
\small
$L_{6a}(2):$ 
$\left(\begin{array}{@{}r@{}r@{}r@{}r@{}r@{}r@{}r@{}r@{}r@{}r@{}r@{}r@{}} 
 2 & 0 & 1 & 1&  0 & 0&  1&  1 & 1 & 1 &\mbox{-}1 &\mbox{-}1\\
 0 & 2 & 0 & 0& \mbox{-}1 & 1&  1&  1 &\mbox{-}1 & 1 & 0 &\mbox{-}1\\
 1 & 0 & 2 & 1&  0 & 0&  1&  1 & 1 & 1 &\mbox{-}1 &\mbox{-}1\\
 1 & 0 & 1 & 2&  0 & 0&  1&  1 & 1 & 1 &\mbox{-}1 &\mbox{-}1\\
 0 &\mbox{-}1 & 0 & 0&  2 &\mbox{-}1&  0& \mbox{-}1 & 1 &\mbox{-}1 & 1 & 1\\
 0 & 1 & 0 & 0& \mbox{-}1 & 2&  0&  1 &\mbox{-}1 & 1 & 0 &\mbox{-}1\\
 1 & 1 & 1 & 1&  0 & 0&  3&  1 & 0 & 1 & 0 &\mbox{-}1\\
 1 & 1 & 1 & 1& \mbox{-}1 & 1&  1&  3 & 0 & 1 &\mbox{-}1 &\mbox{-}1\\
 1 &\mbox{-}1 & 1 & 1&  1 &\mbox{-}1&  0&  0 & 3 & 0 & 0 & 0\\
 1 & 1 & 1 & 1& \mbox{-}1 & 1&  1&  1 & 0 & 3 &\mbox{-}1 &\mbox{-}1\\
\mbox{-}1 & 0 &\mbox{-}1 &\mbox{-}1&  1 & 0&  0& \mbox{-}1 & 0 &\mbox{-}1 & 3 & 1\\
\mbox{-}1 &\mbox{-}1 &\mbox{-}1 &\mbox{-}1&  1 &\mbox{-}1& \mbox{-}1& \mbox{-}1 & 0 &\mbox{-}1 & 1 & 3\end{array} \right)  $,  and 
$L_{6b}(2):$
\small
$\left(\begin{array}{@{}r@{}r@{}r@{}r@{}r@{}r@{}r@{}r@{}r@{}r@{}r@{}r@{}} 
 2&  0 & 0 & 0 & 0 & 0 & 1 & 1 &\mbox{-}1 &\mbox{-}1 & 1 & 1\\
 0&  2 & 0 & 0 & 0 & 0 & 0 & 0 & 0 & 0 & 1 & 1\\
 0&  0 & 2 & 0 & 0 & 0 & 1 &\mbox{-}1 & 0 & 1 &\mbox{-}1 &\mbox{-}1\\
 0&  0 & 0 & 2 & 0 & 0 & 0 & 0 & 1 & 0 & 0 & 0\\
 0&  0 & 0 & 0 & 2 & 0 &\mbox{-}1 & 1 & 0 &\mbox{-}1 & 1 & 1\\
 0&  0 & 0 & 0 & 0 & 2 & 0 & 0 &\mbox{-}1 & 0 & 1 & 1\\
 1&  0 & 1 & 0 &\mbox{-}1 & 0 & 3 & 0 & 0 & 1 &\mbox{-}1 &\mbox{-}1\\
 1&  0 &\mbox{-}1 & 0 & 1 & 0 & 0 & 3 &\mbox{-}1 &\mbox{-}2 & 2 & 2\\
\mbox{-}1&  0 & 0 & 1 & 0 &\mbox{-}1 & 0 &\mbox{-}1 & 3 & 1 &\mbox{-}2 &\mbox{-}2\\
\mbox{-}1&  0 & 1 & 0 &\mbox{-}1 & 0 & 1 &\mbox{-}2 & 1 & 3 &\mbox{-}2 &\mbox{-}2\\
 1&  1 &\mbox{-}1 & 0 & 1 & 1 &\mbox{-}1 & 2 &\mbox{-}2 &\mbox{-}2 & 4 & 3\\
 1&  1 &\mbox{-}1 & 0 & 1 & 1 &\mbox{-}1 & 2 &\mbox{-}2 &\mbox{-}2 & 3 & 4\end{array} \right)$
\\
\normalsize
with automorphism group of order $2^{15} 3^4 $
resp. $2^{21} 3 $.
\item[$L_7(2)$:]   
\small
$\left(\begin{array}{@{}r@{}r@{}r@{}r@{}r@{}r@{}r@{}r@{}r@{}r@{}r@{}r@{}r@{}r@{}} 
  2& 0 & 0 & 0 & 0 & 1 & 1 & 1 & 1& 1&  0 & 0 &-1&  0 \\
  0& 2 & 0 & 0 & 1 & 1 & 1 & 1 & 0& 0&  0 & 0 & 0&  1 \\
  0& 0 & 2 & 0 & 1 & 1 & 0 & 0 &\mbox{-}1& 1&  1 &\mbox{-}1 & 0&  0 \\
  0& 0 & 0 & 2 &\mbox{-}1 & 0 &\mbox{-}1 &\mbox{-}1 & 0& 0& \mbox{-}1 &\mbox{-}1 & 1&  0 \\
  0& 1 & 1 &\mbox{-}1 & 3 & 1 & 1 & 1 & 0& 1&  1 & 0 & 0&  1 \\
  1& 1 & 1 & 0 & 1 & 3 & 1 & 1 & 0& 1&  1 & 0 &\mbox{-}1&  0 \\
  1& 1 & 0 &\mbox{-}1 & 1 & 1 & 3 & 2 & 0& 0&  0 & 0 &\mbox{-}1&  1 \\
  1& 1 & 0 &\mbox{-}1 & 1 & 1 & 2 & 3 & 1& 1&  1 & 1 &\mbox{-}1&  0 \\
  1& 0 &\mbox{-}1 & 0 & 0 & 0 & 0 & 1 & 3& 0&  0 & 1 & 0&  0 \\
  1& 0 & 1 & 0 & 1 & 1 & 0 & 1 & 0& 3&  1 & 0 & 0&  0 \\
  0& 0 & 1 &\mbox{-}1 & 1 & 1 & 0 & 1 & 0& 1&  3 & 0 &\mbox{-}1& \mbox{-}1 \\
  0& 0 &\mbox{-}1 &\mbox{-}1 & 0 & 0 & 0 & 1 & 1& 0&  0 & 3 &\mbox{-}1& \mbox{-}1 \\
 \mbox{-}1& 0 & 0 & 1 & 0 &\mbox{-}1 &\mbox{-}1 &\mbox{-}1 & 0& 0& \mbox{-}1 &\mbox{-}1 & 3&  1 \\
  0& 1 & 0 & 0 & 1 & 0 & 1 & 0 & 0& 0& \mbox{-}1 &\mbox{-}1 & 1&  3 \end{array} \right) $
\normalsize
Automorphism group of order 
2752512.
\item[$L_8(2)$:] 
$L_8(2)$ is the odd version of the Barnes Wall lattice $BW_{16}$ (see 
\cite{database}). It is unique by \cite[Theorem 8]{RaS}.
\item[$L_1(6)$:] 
\small
$\left(\begin{array}{@{}r@{}r@{}r@{}r@{}} 
2 &\mbox{-}1 &\mbox{-}1 & 0\\
\mbox{-}1 & 3 & 0 &\mbox{-}1\\
\mbox{-}1 & 0 & 3 &\mbox{-}1\\
 0 &\mbox{-}1 &\mbox{-}1 & 4\end{array} \right) $.
\normalsize
Automorphism group $C_2^4$.
\item[$L_2(6)$:]
\small
$\left(\begin{array}{@{}r@{}r@{}r@{}r@{}r@{}r@{}r@{}r@{}} 
 3 & 1 & 0 & 0 & 0 & 1 & 0 &\mbox{-}1\\
 1 & 3 & 0 & 0 &\mbox{-}1 & 0 &\mbox{-}1 & 0\\
 0 & 0 & 3 & 1 & 0 &\mbox{-}1 & 0 &\mbox{-}1\\
 0 & 0 & 1 & 3 & 1 & 0 &\mbox{-}1 & 0\\
 0 &\mbox{-}1 & 0 & 1 & 3 & 1 & 0 & 0\\
 1 & 0 &\mbox{-}1 & 0 & 1 & 3 & 0 & 0\\
 0 &\mbox{-}1 & 0 &\mbox{-}1 & 0 & 0 & 3 &\mbox{-}1\\
\mbox{-}1 & 0 &\mbox{-}1 & 0 & 0 & 0 &\mbox{-}1 & 3\end{array} \right)$
\normalsize
Automorphism group $SL_2(3).2^2$ of order $96$.
\item[$L_1(14)$:]
Gram matrix $\left( \begin{array}{cc} 3 & 1 \\ 1 & 5 \end{array} \right) \perp
\left( \begin{array}{cc} 3 & 1 \\ 1 & 5 \end{array} \right) $.
Automorphism group $D_8$.
\end{itemize}

\end{document}